\documentclass[11pt]{amsart}

\newif\ifdraft\draftfalse

\usepackage{amssymb}
\usepackage{amsthm}
\usepackage{eucal}
\usepackage[english]{babel}
\usepackage{nicefrac}
\usepackage{times}
\usepackage{latexsym,amscd,amsmath,amsthm,amssymb}

\makeatletter
\def\@begintheorem#1#2[#3]{%
    \def\naam{#1}
  \deferred@thm@head{\the\thm@headfont \thm@indent
    \@ifempty{#1}{\let\thmname\@gobble}{\let\thmname\@iden}%
    \@ifempty{#2}{\let\thmnumber\@gobble}{\let\thmnumber\@iden}%
    \@ifempty{#3}{\let\thmnote\@gobble}{\let\thmnote\@iden}%
    \thm@swap\swappedhead\thmhead{#1}{#2}{#3}%
    \the\thm@headpunct
    \thmheadnl 
    \hskip\thm@headsep
  }%
  \ignorespaces}
\makeatother

\newcommand{\kantlijndraft}[1]{\ifdraft\hspace{-\lastskip}%
\vadjust{\vspace{-1mm}\smash{\llap{{\tt #1}\hspace{8mm}}}\vspace{1mm}}\fi}

\def\voegToe#1#2#3{\immediate\write1{\string\newlabel{#1}{{#2}{#3}}}}

\newcommand{\thlabel}[1]{\voegToe{#1}{\naam\noexpand~\thetheorem}{\thepage}\kantlijndraft{#1}}

\makeatletter
\renewcommand{\label}[1]{\voegToe{#1}{\@currentlabel}{\thepage}\kantlijndraft{#1}}
\makeatother

\newtheorem{theorem}{Theorem}[section]

\newtheorem{corollary}[theorem]{Corollary}
\newtheorem{question}[theorem]{Question}

\newtheorem{proposition}[theorem]{Proposition}
\theoremstyle{definition}
\newtheorem{example}[theorem]{Example}
\newtheorem{definition}[theorem]{Definition}
\theoremstyle{remark}

\numberwithin{equation}{section}

\newtheorem{claim2}{\sc Claim}



\newcommand{\sse}{\subseteq}						
\newcommand{\minus}{\backslash}						
\newcommand{\Un}{\bigcup}							
\newcommand{\un}{\cup}								
\newcommand{\Meet}{\bigcap}							
\newcommand{\meet}{\cap}							
\newcommand{\es}{\varnothing}						

\newcommand{\cl}[1]{\ensuremath{\overline{#1}}}

\newcommand{\scr}[1]{\ensuremath{\mathcal{#1}}}

\def\sapirovskii{{\v{S}}apirovski{\u\i}}

\def\juhasz{Juh{\'a}sz}

\begin{document}

\title{A bound for the density of any Hausdorff space}

\author{Nathan Carlson}\address{Department of Mathematics, California Lutheran University, 60 W. Olsen Rd, MC 3750, 
Thousand Oaks, CA 91360 USA}
\email{ncarlson@callutheran.edu}

\begin{abstract}
We show, in a certain specific sense, that both the density and the cardinality of a Hausdorff space are related to the ``degree" to which the space is nonregular. It was shown by \sapirovskii~that $d(X)\leq\pi\chi(X)^{c(X)}$ for a regular space $X$ and the author observed this holds if the space is only quasiregular. We generalize this result to the class of all Hausdorff spaces by introducing the nonquasiregularity degree $nq(X)$, which is countable when $X$ is quasiregular, and showing $d(X)\leq\pi\chi(X)^{c(X)nq(X)}$ for any Hausdorff space $X$. This demonstrates that the degree to which a space is nonquasiregular has a fundamental and direct connection to its density and, ultimately, its cardinality. Importantly, if $X$ is Hausdorff then $nq(X)$ is ``small" in the sense that $nq(X)\leq\psi_c(X)$. This results in a unified proof of both \sapirovskii's density bound for regular spaces and Sun’s bound $\pi\chi(X)^{c(X)\psi_c(X)}$ for the cardinality of a Hausdorff space $X$. A consequence is an improved bound for the cardinality of a Hausdorff space.
\end{abstract}

\subjclass[2020]{54A25, 54D10.}

\keywords{cardinality bounds, cardinal functions}

\maketitle

\section{Introduction.}

The \emph{density} of a topological space $X$, denoted by $d(X)$, is the least cardinality of a dense subset of $X$. In 1974 \sapirovskii~showed that the density of a regular space $X$ is bounded above by $\pi\chi(X)^{c(X)}$ \cite{Sap1974}. It was observed by the author in \cite{C2007a} that this holds if the space has a weaker property called quasiregular (Definition~\ref{qr}) and not necessarily even Hausdorff. Additional proofs of this result were given in~\cite{Cha77} and~\cite{Juh80}, all of which are sophisticated arguments. In 1988 Sun showed $|X|\leq\pi\chi(X)^{c(X)\psi_c(X)}$ for any Hausdorff space X, again using a complicated closing-off argument, but quite different than the proof of \sapirovskii’s density bound. 

We establish a bound for the density of any space in a certain class that includes Hausdorff spaces and quasiregular spaces. This results in a unified proof of \sapirovskii's and Sun’s theorems, gives an improved bound for the cardinality of a Hausdorff space, and answers Question 3.3 in~\cite{C2023b}. To achieve this we introduce the nonquasiregular degree $nq(X)$ (Definition~\ref{nq}) which gives a ``measure” of how nonquasiregular a space is. A quasiregular space X has countable $nq(X)$. In addition, and importantly, if $X$ is Hausdorff then $nq(X)\leq\psi_c(X)$. This shows that the nonquasiregular degree generalizes both quasiregularity in spaces for which $nq(X)$ is defined and the closed pseudocharacter in Hausdorff spaces. 

Our main result is that $d(X)\leq\pi\chi(X)^{c(X)nq(X)}$ for any space $X$ for which $nq(X)$ is defined, including Hausdorff spaces and quasiregular spaces. That is, such spaces have a dense subset $D$ such that $|D|\leq \pi\chi(X)^{c(X)nq(X)}$. This shows, in a sense, that $nq(X)$ is the ``missing piece of the puzzle" in \sapirovskii's density bound. A corollary is that $|X|\leq\pi\chi(X)^{c(X)nq(X)w\psi_c(X)}$  for any Hausdorff space $X$, where $w\psi_c(X)$ was introduced in~\cite{C2023b} and satisfies $w\psi_c(X)\leq\psi_c(X)$ if $X$ is Hausdorff. (See Definition~\ref{wpsi}). This improves Sun’s theorem and shows that the degree to which a space is nonquasiregular plays an important role in the cardinality of a Hausdorff space. 

It is also shown that if $X$ is Hausdorff then $nq(X)\leq L(X)$. This gives, in addition to quasiregular spaces, a new class of spaces for which~\sapirovskii's inequality $d(X)\leq\pi\chi(X)^{c(X)}$ holds: Lindel\"of Hausdorff spaces.

We give two examples of spaces $X$ where $nq(X)$ is required in our density bound; that is, $d(X)>\pi\chi(X)^{c(X)}$. One of these examples is even Urysohn. We also give several other examples.

No global assumptions are made on any separation axiom on a space.  For definitions of cardinal functions not defined here, see \juhasz~\cite{Juh80}.

\section{A bound for the density.}

\begin{definition}\label{qr}
A space $X$ is \emph{quasiregular} if every nonempty open set $U$ of $X$ contains a nonempty open set $V$ such that $\cl{V}\sse U$.
\end{definition}

Note that any space with a dense set of isolated points is quasiregular. We introduce the nonquasiregularity degree $nq(X)$ for a space $X$.

\begin{definition}\label{nq}
Let $X$ be a space. The \emph{nonquasiregularity degree} $nq(X)$ is defined as the least infinite cardinal $\kappa$ such that every nonempty open set $U$ contains a nonempty $G_\kappa$-set $G=\Meet\scr{V}$ such that $\scr{V}$ is an open family, $|\scr{V}|\leq\kappa$, and $\Meet_{V\in\scr{V}}\cl{V}\sse U$.
\end{definition}

Observe that $nq(X)$ is defined when $X$ is Hausdorff and when $X$ is quasiregular. In the latter case we have $nq(X)=\aleph_0$. However, there are spaces for which $nq(X)$ is not defined. Example~\ref{notdefined} gives an example of compact, $T_1$ space for which $nq(X)$ is not defined. 

\begin{definition} A space $X$ is an $nq$-\emph{space} if $nq(X)$ is defined.
\end{definition}

Every Hausdorff space is an $nq$-space and every quasiregular space is an $nq$-space. 

The next two definitions were given in \cite{C2023b}.


\begin{definition}\label{wpsi}
Let $X$ be a Hausdorff space and let $x\in X$. A collection of open sets $\scr{V}$ is a \emph{weak closed pseudobase at }$x$ if $\{x\}=\Meet_{V\in\scr{V}}\cl{V}$. (Note that it is not necessarily the case that $x\in V$ for any $V\in\scr{V}$). We define $w\psi_c(x,X)$ to be the least infinite cardinal $\kappa$ such that $x$ has a weak closed pseudobase of cardinality $\kappa$. The \emph{weak closed pseudocharacter} $w\psi_c(X)$ is defined as $w\psi_c(X)=\sup\{w\psi_c(x,X):x\in X\}$.
\end{definition}

\begin{definition}\label{dpsi}
Let $X$ be a Hausdorff space. The \emph{dense closed pseudocharacter} $d\psi_c(X)$ of $X$ is defined as the least infinite cardinal $\kappa$ such that $X$ has a dense set $D$ such that $\psi_c(d,X)\leq\kappa$ for every $d\in D$.
\end{definition}

Observe that $nq(X)\leq d\psi_c(X)\leq\psi_c(X)$ if $X$ is Hausdorff. We prove now the central result in this paper.

\begin{theorem}
Let $X$ be an $nq$-space. Then $d(X)\leq\pi\chi(X)^{c(X)nq(X)}$.
\end{theorem}

\begin{proof}
Let $\lambda=\pi\chi(X)$ and $\kappa=c(X)nq(X)$. For all $x\in X$ let $\scr{B}_x$ be a local $\pi$-base for $x$ such that $|\scr{B}_x|\leq\lambda$. By transfinite induction we construct a non-decreasing chain of $\{A_\alpha:\alpha<\kappa^+\}$ of subsets of $X$ and a sequence of open collections $\{\scr{B}_\alpha:\alpha<\kappa^+\}$ such that the following properties hold for all $\alpha<\kappa^+$:

\begin{enumerate}
\item $|A_\alpha|\leq\lambda^\kappa$,
\item $|\scr{B}_\alpha|\leq\lambda^\kappa$, and 
\item if $\scr{U}=\{\scr{U}_\gamma:\gamma<\kappa\}\in\left[[\scr{B}_\alpha]^{\leq\kappa}\right]^{\leq\kappa}$ and $X\minus\Un_{\gamma<\kappa}\cl{\Un\scr{U}_\gamma}\neq\es$, then $A_\alpha\minus\Un_{\gamma<\kappa}\cl{\Un\scr{U}_\gamma}\neq\es$. 
\end{enumerate}

Pick $p\in X$. Let $A_0=\{p\}$ and $\scr{B}_0=\scr{B}_p$. Let $0<\alpha<\kappa^+$ and assume that $\{A_\beta:\beta<\alpha\}$ have been constructed. Define $\scr{B}_\alpha=\Un\{\scr{B}_x:x\in\Un_{\beta<\alpha}A_\beta\}$. Then $|\scr{B}_\alpha|\leq\lambda\cdot\lambda^\kappa\cdot\kappa=\lambda^\kappa$. For each $\scr{U}=\{\scr{U}_\gamma:\gamma<\kappa\}\in\left[[\scr{B}_\alpha]^{\leq\kappa}\right]^{\leq\kappa}$ such that $X\minus\Un_{\gamma<\kappa}\cl{\Un\scr{U}_\gamma}\neq\es$, pick $x_\scr{U}\in X\minus\Un_{\gamma<\kappa}\cl{\Un\scr{U}_\gamma}$. Define 
$$A_\alpha=\Un_{\beta<\alpha}A_\beta\un\left\{x_\scr{U}:\scr{U}=\{\scr{U}_\gamma:\gamma<\kappa\}\in\left[[\scr{B}_\alpha]^{\leq\kappa}\right]^{\leq\kappa}\textup{ such that }X\minus\Un_{\gamma<\kappa}\cl{\Un\scr{U}_\gamma}\neq\es\right\}.$$
As $\left|\Un_{\beta<\alpha}A_\beta\right|\leq \lambda^\kappa\cdot\kappa=\lambda^\kappa$ and $\left|\left[[\scr{B}_\alpha]^{\leq\kappa}\right]^{\leq\kappa}\right|\leq \left((\lambda^\kappa)^\kappa\right)^\kappa=\lambda^\kappa$, we see that $|A_\alpha|\leq\lambda^\kappa$. By the way we have constructed $A_\alpha$ we see that (3) is satisfied.

Let $A=\Un_{\alpha<\kappa^+}A_\alpha$. Then $|A|\leq\kappa^+\lambda^\kappa=\lambda^\kappa$. We show that $A$ is dense in $X$. Suppose by way of contradiction that there exists a nonempty open set $U$ such that $U\meet A=\es$. As $nq(X)\leq\kappa$, there exists a nonempty $G_\kappa$-set $G=\Meet\scr{V}$ such that $\scr{V}=\{V_\alpha:\alpha<\kappa\}$ is an open family and $\Meet_{\alpha<\kappa}\cl{V_\alpha}\sse U$.

For every $\alpha<\kappa$, let $W_\alpha=X\minus\cl{V_\alpha}$. Observe that $\cl{W_\alpha}=cl(X\minus\cl{V_\alpha})=X\minus int(\cl{V_\alpha})\sse X\minus V_\alpha\sse X\minus G$ and so $G\meet\cl{W_\alpha}=\es$. Also, $A\sse X\minus U\sse X\minus\Meet_{\alpha<\kappa}\cl{V_\alpha}=\Un_{\alpha<\kappa}X\minus\cl{V_\alpha}=\Un_{\alpha<\kappa}W_\alpha$ and so $A\sse \Un_{\alpha<\kappa}W_\alpha$. For each $\alpha<\kappa$, define $\scr{S}_\alpha=\{B\in\scr{B}_x: x\in W_\alpha\meet A, B\sse W_\alpha\}$. Note $\Un\scr{S}_\alpha\sse W_\alpha$.

We show $A\meet W_\alpha\sse\cl{\Un\scr{S}_\alpha}$ for all $\alpha<\kappa$. Let $x\in A\meet W_\alpha$ and let $T$ be an open set containing $x$. There exists $B\in\scr{B}_x$ such that $B\sse W_\alpha\meet T$. Then $B\in\scr{S}_\alpha$ and $\es\neq B\sse T\meet\Un\scr{S}_\alpha$. This shows $A\meet W_\alpha\sse\cl{\Un\scr{S}_\alpha}$ for all $\alpha<\kappa$.

As $c(X)\leq\kappa$, using maximal pairwise disjoint families, for each $\alpha<\kappa$ there exists $\scr{U}_\alpha\sse\scr{S}_\alpha$ such that $\cl{\Un\scr{U}_\alpha}=\cl{\Un\scr{S}_\alpha}$ and $|\scr{U}_\alpha|\leq\kappa$. Now, for each $\alpha<\kappa$ note that $\cl{\Un\scr{U}_\alpha}=\cl{\Un\scr{S}_\alpha}\sse\cl{W_\alpha}\sse X\minus G$. Therefore, $\es\neq G\sse X\minus\Un_{\alpha<\kappa}\cl{\Un\scr{U}_\alpha}$. Since $\left|\Un_{\alpha<\kappa}\scr{U}_\alpha\right|\leq\kappa\cdot\kappa=\kappa<\kappa^+$, there exists $\delta<\kappa^+$ such that $\scr{U}=\{\scr{U}_\alpha:\alpha<\kappa\}\in\left[[\scr{B}_\delta]^{\leq\kappa}\right]^{\leq\kappa}$.

By (3) above, we have that $x_\scr{U}\in A_{\delta+1}\minus\Un_{\alpha<\kappa}\cl{\Un\scr{U}_\alpha}\sse A\minus\Un_{\alpha<\kappa}\cl{\Un\scr{U}_\alpha}$. This contradicts the fact that $A\sse\Un_{\alpha<\kappa}\cl{\Un\scr{S}_\alpha}=\Un_{\alpha<\kappa}\cl{\Un\scr{U}_\alpha}$. Therefore $A$ is dense in $X$ and
$d(X)\leq |A|\leq\lambda^\kappa=\pi\chi(X)^{c(X)nq(X)}$.
\end{proof}

As every Hausdorff space is an $nq$-space and every quasiregular space is an $nq$-space, we have the following corollary. This provides a unified proof that $d(X)\leq\pi\chi(X)^{c(X)}$ for quasiregular spaces and that $d(X)\leq\pi\chi(X)^{c(X)d\psi_c(X)}$ for Hausdorff spaces \cite{C2023b}, answering Question 3.3 in \cite{C2023b}. 

\begin{corollary}\label{dense}
Let $X$ be Hausdorff or quasiregular. Then $d(X)\leq\pi\chi(X)^{c(X)nq(X)}$.
\end{corollary}

\begin{corollary}
Let $X$ be a space.
\begin{enumerate}
\item (\sapirovskii~in the regular case~\cite{Sap1974} , C. for quasiregular~\cite{C2007a}) If $X$ is quasiregular then $d(X)\leq\pi\chi(X)^{c(X)}$, and
\item (C. \cite{C2023b}) If $X$ is Hausdorff then $d(X)\leq\pi\chi(X)^{c(X)d\psi_c(X)}$.
\end{enumerate}
\end{corollary}

Recall a set $R$ in a space $X$ is \emph{regular open} if $R=int\cl R$. Let $RO(X)$ denote the set of regular open sets of a space $X$.

\begin{corollary}\label{ro}
Let $X$ be an $nq$-space. Then $|RO(X)|\leq\pi\chi(X)^{c(X)nq(X)}$.
\end{corollary}

\begin{proof}
It is well-known that $|RO(X)|\leq\pi w(X)^{c(X)}$ for any space $X$. We have,
\begin{align}
|RO(X)|\leq\pi w(X)^{c(X)}=(d(X)\pi\chi(X))^{c(X)}&\leq\left(\pi\chi(X)^{c(X)nq(X)}\cdot\pi\chi(X)\right)^{c(X)}\notag\\
&=\pi\chi(X)^{c(X)nq(X)}.\notag
\end{align}
\end{proof}

\begin{corollary}\label{Hausdorffbound}
For any Hausdorff space $X$, $|X|\leq\pi\chi(X)^{c(X)nq(X)w\psi_c(X)}$.
\end{corollary}

\begin{proof}
Apply the above Corollary~\ref{ro} and the fact that $|X|\leq|RO(X)|^{w\psi_c(X)}$, as shown in \cite{C2023b}.
\end{proof}

\begin{corollary}
Let $X$ be a space.
\begin{enumerate}
\item (C.~\cite{C2023b}) If $X$ Hausdorff then $|X|\leq\pi\chi(X)^{c(X)d\psi_c(X)w\psi_c(X)}$.
\item (Sun~\cite{Sun88}) If $X$ is Hausdorff then $|X|\leq\pi\chi(X)^{c(X)\psi_c(X)}$.
\item (C.~\cite{C2023b}) If $X$ is quasiregular and Hausdorff then $|X|\leq\pi\chi(X)^{c(X)w\psi_c(X)}$.
\item (\sapirovskii~\cite{Sap1974}) If $X$ is regular and Hausdorff then $|X|\leq\pi\chi(X)^{c(X)\psi(X)}$.
\end{enumerate}
\end{corollary}

We conclude this section by showing that $nq(X)\leq L(X)$ if $X$ is Hausdorff.

\begin{proposition}\label{Lindelof}
If $X$ is Hausdorff then $nq(X)\leq L(X)$.
\end{proposition}

\begin{proof}
Let $\kappa=L(X)$. Let $U$ be a nonempty open set and let $p\in U$. As $X$ is Hausdorff, for every $x\in X\minus U$ there exists an open set $V_x$ containing $x$ such that $p\in X\minus\cl{V_x}$. As $X\minus U$ is closed, there exists $A\in[X\minus U]^{\leq\kappa}$ such that $X\minus U\sse\Un_{x\in A}V_x$. Thus, 
$$p\in\Meet_{x\in A}X\minus\cl{V_x}\sse\Meet_{x\in A}cl(X\minus\cl{V_x})=\Meet_{x\in A}X\minus int\cl{V_x}\sse\Meet_{x\in A}X\minus V_x\sse  U.$$
This shows every nonempty open set contains a nonempty $G_\kappa$-set $G=\Meet\scr{V}$ such that $\Meet_{V\in\scr{V}}\cl{V}\sse U$. Therefore, $nq(X)\leq\kappa=L(X)$.
\end{proof}

By Corollary~\ref{dense} and Proposition~\ref{Lindelof}, we have a new class of spaces for which~\sapirovskii's inequality holds: Lindel\"of Hausdorff spaces. This appears to be new in the literature.

\begin{corollary} If $X$ is Lindel\"of and Hausdorff then $d(X)\leq\pi\chi(X)^{c(X)}$.
\end{corollary}

\section{Examples.}

In this section we give several examples. Examples~\ref{exampleA} and~\ref{exampleB} are examples of spaces $X$ where $d(X)>\pi\chi(X)^{c(X)}$, demonstrating that $nq(X)$ is necessary in Theorem~\ref{dense}. Example~\ref{exampleB} is even Urysohn.

\begin{example}\label{exampleA}
In~\cite{Cha77} Charlesworth attributes this example to van Douwen. Let $\kappa$ be an infinite cardinal and let $X$ be the product of $2^\kappa$ copies of the two point discrete space and let the topology on $X$ be defined as follows. Let $\leq$ be a well-ordering on $X$. For every $A\sse X$ and $x\in X$, let $A_x=\{y\in A: x\leq y\}$. Basic open sets are of the form $U_x$, where $x\in X$ and $U$ is open in the usual product topology on $X$. It can be seen that $X$ is Hausdorff, $\chi(X)=2^\kappa$ and $d(X)>cf(2^{2^\kappa})>2^\kappa$. A proof that $X$ is c.c.c. is in~\cite{Cha77}. This shows $d(X)>2^\kappa=\chi(X)^{c(X)}\geq\pi\chi(X)^{c(X)}$. 

However, as $d(Y)\leq\pi\chi(Y)^{c(Y)nq(Y)}$ for any space $X$, we have that $d(X)\leq\pi\chi(X)^{c(X)nq(X)}\leq 2^{\kappa\cdot nq(X)}$. This implies $nq(X)\geq\kappa^+$ and that $nq(X)$ is necessary in our density bound.
\end{example}

The next example requires Proposition 4.1 from~\cite{CR2008}. Recall that $X_s$ represents the \emph{semiregularization} of a space $X$. (See~\cite{PW}).

\begin{proposition}\label{prop}
Fix an infinite cardinal number $\kappa$ and let $Z$ be any space of size $\kappa+$ such that $\pi\chi(Z)\leq\kappa$. Furthermore, assume that $|U|=\kappa^+$ whenever $U$ is a nonempty open subset of $Z$. Then there is a space $X$ such that $X_s=Z_s$, $\pi\chi(X)=\kappa$, and $d(X)=\kappa^+$.
\end{proposition}

\begin{example}\label{exampleB}
We give another example of a space $X$ in which $nq(X)$ is necessary in Theorem~\ref{dense}. Assume $\mathfrak{c}^+=2^\mathfrak{c}$ and let $Z=\mathbb{I}^{\mathfrak{c}}$. Apply Proposition~\ref{prop} to obtain a space $X$ that $\pi\chi(X)=\mathfrak{c}$ and $d(X)=\mathfrak{c}^+$. Furthermore, since $c(X)=c(X_s)=c(Z)$ and $Z$ is c.c.c, then $X$ is c.c.c. In addition, $X$ is Urysohn as $Z$ is Urysohn and $X_s=Z$. We have $d(X)=\mathfrak{c}^+=2^\mathfrak{c}\leq\pi\chi(X)^{c(X)nq(X)}=\mathfrak{c}^{\omega\cdot nq(X)}=2^{nq(X)}$. This implies that $nq(X)\geq\mathfrak{c}$ and that $nq(X)$ is necessary in our density bound.
\end{example}

We give an example of a compact, $T_1$ space that is not an $nq$-space. This example was independently noticed by Angelo Bella and the author.

\begin{example}\label{notdefined}
Consider the cofinite topology on an infinite set $X$. This space is well-known to be compact, $T_1$, and not Hausdorff. Let $F$ be a nonempty finite set and consider the open set $X\minus F$. Let $\kappa$ be any cardinal and consider any $G_\kappa$ set $\Meet\scr{V}\sse X\minus F$ where $\scr{V}$ is an open family and $|\scr{V}|\leq\kappa$. Now as the space is hyperconnected, we have that every nonempty open set is dense. This means $\cl{V}=X$ for every $V\in\scr{V}$. Therefore $\Meet_{V\in\scr{V}}\cl{V}=X$ and $\Meet_{V\in\scr{V}}\cl{V}$ is not contained in $X\minus F$. This shows $nq(X)$ is not defined.
\end{example}

The next example is an example of a nonquasiregular space $X$ for which $nq(X)$ is countable. This example was communicated by Jack Porter.

\begin{example}\label{porter}
Let $X=\{(x,y):y\geq 0, x, y\in\mathbb{Q}\}$ with the well-known ``irrational slope topology". (See~\cite{SS1978}). This is also called ``Bing's Tripod Space." This space is easily seen to be nonquasiregular. In addition it is Hausdorff and first countable, and thus $nq(X)\leq\psi_c(X)=\aleph_0$. 
\end{example}

Finally, we observe that any set with two or more points with the indiscrete topology is regular, hence quasiregular with countable nonquasiregularity degree. This shows there are spaces where $nq(X)$ is defined that are not even $T_0$.

In light of Theorem~\ref{Hausdorffbound}, we ask the following.

\begin{question}
Is there a Hausdorff space $X$ such that $\pi\chi(X)^{c(X)nq(X)w\psi_c(X)}<\pi\chi(X)^{c(X)\psi_c(X)}$?
\end{question}

We also ask if we can distinguish $nq(X)$ and $d\psi_c(X)$ in Hausdorff spaces.

\begin{question} Is there a Hausdorff space for which $nq(X)<d\psi_c(X)$?
\end{question}

We know $nq(X)$ is defined for Hausdorff spaces and quasiregular spaces. Related to this we ask the following question.

\begin{question}
Can we characterize the $nq$-spaces?
\end{question}

\textbf{Acknowledgement.} The author would like to thank Jack Porter and Angelo Bella for a careful reading of this manuscript and helpful discussion.

\end{document}